
\def\khetzi{ {{1} \over {2}}}

\baselineskip=14pt
\parskip=10pt

\font\eightrm=cmr8 
\font\eighttt=cmtt8
\magnification=\magstephalf

\def\1{{\overline{1}}}
\def\2{{\overline{2}}}
\def\Tilde{\char126\relax}
\parindent=0pt
\overfullrule=0in
\def\frac#1#2{{#1 \over #2}}
\bf
\centerline
{
Searching for Strange Hypergeometric Identities By Sheer Brute Force
}
\rm
\bigskip
\centerline 
{ 
{\it Moa APAGODU}
and {\it Doron ZEILBERGER}\footnote{$^1$}
{
{\eightrm
Supported in part by the USA National Science Foundation.
}
}
}
\bigskip

{\bf Important Note:} This article is accompanied by the Maple
package {\tt BruteTwoFone} available from the
webpage of this article \hfill\break
{\eighttt http://www.math.rutgers.edu/\Tilde zeilberg/mamarim/mamarimhtml/sheer.html} \quad , \hfill\break
where one can also find sample input and output.

{\bf Preface}

The classical {\it hypergeometric} series
$$
F(a\, , \, b \, , \, c \, , \,x) \, = \,
\sum_{k=0}^{\infty}
{ {(a)_k (b)_k}  \over {k! (c)_k } } x^k \quad,
$$
(where $(z)_k:=z(z+1)(z+2) \cdots (z+k-1)$),
that nowadays is more commonly denoted by
$$
{}_2 F_1 \left ( {{a,b} \atop {c}} \, ; \, x \right ) \quad,
$$
has a long and distinguished history, going back to Lehonard Euler
and Carl Friedrich Gauss. It was also one of Ramanujan's favorites.
Under the guise of {\bf binomial
coefficient sums} it goes even further back, to Chu, in
his 1303 combinatorics treatise, that summarized a body of
knowledge that probably goes yet further back.

The hypergeometric function, and its generalized counterparts,
enjoy several {\it exact evaluations}, for some choices of
the parameters, in terms of the Gamma function. The classical
identities of Chu-Vandermonde, Gauss, Kummer, Euler, 
Pfaff-Saalsch\"utz, Dixon, Dougall, and others, can be looked
up in the {\it classic} classic text of Bailey[B], and the 
{\it modern} classic text of Andrews, Askey, and Roy [AAR].

For example, when $x=1$, Gauss found the 
$3$-parameter exact evaluation:
$$
F(a\, ,\, b \, , \, c \, , \, 1) \, =  \,
{ {\Gamma(c) \Gamma(c-a-b)} \over {\Gamma(c-a) \Gamma(c-b)}} \quad .
\eqno(Gauss)
$$
(When $a$ is a negative integer, $a=-n$, then this goes back
to Chu's 1303 identity, rediscovered by Vandermonde).

Next comes Kummer's {\it two-parameter } exact evaluation at
$x=-1$
$$
F(a\, ,\, b \, , \, 1+a-b \, , \, -1)= 
{{\Gamma(1+a-b) \Gamma(1+a/2)} \over {\Gamma(1+a) \Gamma(1+a/2-b)}} \quad ,
\eqno(Kummer)
$$
and two other ones, at $x=\khetzi$, due to Gauss (see [B] or [AAR]).
In addition, Gosper conjectured,
and Gessel and Stanton[GS] proved, several `strange'
{\it one-parameter} evaluations at other values of $x$.

The Hypergeometric function also enjoys several {\it transformation
formulas}, both rational and quadratic, due to Euler, and Pfaff
(see [AAR], Theorem 2.2.5, Corollary 2.3.3 and Theorem 3.1.3),
so any  exact evaluation implies quite a few other ones, equivalent
to it via these transformations, and iterations thereof.
[See procedures {\tt Buddies21C} and {\tt QuadBuddies21}, in our Maple
package {\tt BruteTwoFone}. ]

In [E], a systematic search for {\it all} such strange 
identities, up to a certain ``complexity''  was attempted. 
It was done by implementing the method of [Z], that
used Wilf-Zeilberger theory and Zeilberger's
algorithm as simplified in [MZ]. The drawback of that method, however,
was that it only searched for those identities for which
the Zeilberger algorithm outputs a {\it first-order} recurrence.
While rare, there are cases where the Zeilberger algorithm outputs
a higher-order recurrence, yet still is evaluable in closed-form.
This is because this algorithm is {\bf not} guaranteed to output
the minimal-order recurrence, although it usually does.

This gave us the idea to systematically search for such 
closed-form evaluations by {\it sheer brute force}.
Surprisingly, it lead us to the discovery of two {\it new}
infinite families of ``closed-form'' evaluations,
that we will describe later.

When $a$ or $b$ happen to be a negative integer,$a=-n$, say, then
the {\it infinite} series {\bf terminates}, for the
sum then only has $n+1$ terms, and one does not have to worry about
convergence. In most cases, one can easily pass from
the terminating case to the non-terminating case by 
Carlson's theorem ([AAR], p. 108; [B], p. 39). In this article
we will only consider such terminating series.

{\bf The Haystack}

In this {\it etude} in {\bf Experimental Mathematics}, the haystack
consists of 
$$
F(-an, bn+b_0, cn+c_0, x) \quad ,
\eqno(2F1)
$$
where $a$ is a positive integer, $b$ and $c$ are integers, while 
$b_0,c_0,x$ are complex numbers.

{\bf The Needles}

The needles are those  $(2F1)$'s that
are evaluable in terms of the Gamma function,
or more precisely, those for which the sequence
$$
u_n:=F(-an, bn+b_0, cn+c_0, x) \quad ,
$$
is a {\it hypergeometric} sequence which means that
$$
r_n:={{u_{n+1}} \over {u_n}} \quad ,
$$
is a {\it rational function} of $n$, i.e.
$$
r_n={{P(n)} \over {Q(n)}} \quad,
$$
where $P(n)$ and $Q(n)$ are {\it polynomials} in $n$.
It is easy to see, by looking at the asymptotics, that
the degrees of $P$ and $Q$ must be the same.

{\bf How to Test Whether Something is a Needle or just a 
Boring piece of Hay?}

Even with a {\it specific} choice of $a,b,c,b_0,c_0,x$, there
is no way, {\it a priori}, to rule out (conclusively) whether the resulting
sequence, $u_n$, is hypergeometric, since, who knows?, the
degree $d$ could be a zillion. But if we restrict the
search for some fixed (not too big!) $d$, then it is
very easy (with computers, of course), to decide whether
the studied sequence $u_n$ is hypergeometric with the
degree of both top and bottom of the $r_n$ being $\leq d$.

Indeed, write $r_n=P(n)/Q(n)$ {\it generically}, 
$$
r_n={ {\sum_{i=0}^{d} p_i n^i } \over
{\sum_{i=0}^{d} q_i n^i } } \quad ,
$$
in terms of the $2d+2$ {\bf undetermined coefficients}
$p_0,p_1, \dots p_d$, $q_0, q_1, \dots , q_d$, and plug-in
$n=0,1,2, \dots, 2d+6$, say. You would get $2d+7$ equations
for the $2d+2$ unknowns, and short of a miracle, they would
not be solvable. If there are solvable, then with probability
$1-10^{-10000}$, $r_n$ is indeed the conjectured rational
function, and if you want to have it true with probability $1$, then
all you need to do is find the corresponding $u_n$
(that solves $u_{n+1}/u_n=r_n$), and then use Zeilberger's
algorithm to prove the conjectured identity rigorously.

[The above is done in procedure {\tt NakhD} in our Maple
package {\tt BruteTwoFone}. For example, to guess the
Chu-Vandermonde identity, type
{\tt NakhD([[-n,a],[c],1],n,1);}, and you would get \hfill\break
$(c+n-a)/(c+n)$ ]

{\bf Alas the Haystack is infinite}

Even with a specific $a$ and a specific degree $d$,
there are {\it infinitely} many things to try.
However, if we leave $b,b_0,c,c_0,x$ {\it symbolic}, then
for any specific, numeric $n_0$, $u_{n_0}$ will no longer
be a mere number, but a certain
{\it rational function} of  $(b,b_0,c,c_0,x)$.

We are looking for those choices of the parameters 
$(b,b_0,c,c_0,x)$, for which
the linear equations in the coefficients of $r_n$, namely
$p_0, \dots , p_d$, $q_0, \dots , q_d$, are solvable.
This means, using linear algebra, that certain determinants
vanish. We can add as many conditions as we want, by
finding the determinants corresponding to
the set of equation $P(i)/Q(i)=r_i$, for $i=1, \dots, M$,
for $M$ big enough to have more equations than unknowns.
Then using the Buchberger algorithm, we can solve them, and
get {\it all} the choices of $b,b_0,c,c_0,x$
(including infinite families, for example, we should get
$x=1$ to account for the Chu-Vandermonde identity).

{\bf Alas our Computers are not Big Enough}

Unfortunately, the above scheme is not feasible, since
the equations are {\it sooo} huge, and Buchberger's algorithm
is {\it sooo} slow. So we have to compromise.
We {\it fix} $b,b_0,c,c_0$, and search for the lucky $x$.
Now we only need two determinants, both being
certain polynomials in the
single variable $x$, and simply take their greatest common
divisor ({\tt gcd} in Maple), to get those $x$ that 
(have the potential, and probably will) yield
closed-form evaluations (with the given $d$) for those
{\it fixed} $a,b,c,b_0,c_0$. 

[The above is done in procedure {\tt NakhDx} in our Maple
package {\tt BruteTwoFone}. For example, to guess Theorem 12
of [E], type

{\tt NakhDx([-n,-3*n-1],[-2*n],n,10);}, 

and you would  get
$x=(1 \pm \sqrt{3} i)/2$, for the two choices for $x$ that would make \hfill\break
$F(-n, -3n-1, -2n, x)$ hypergeometric of degree $\leq 10$. ]

{\bf The Big Five-Fold Do-Loop}

So we decide beforehand on an $M$, and a denominator $D$, and try
{\it all} $F(-an,bn+b_0,cn+c_0,x)$ ($x$ yet-to-be-determined)
with integers $a,b,c$ such that  $1 \leq a \leq M$, and $-M \leq b ,c \leq M$, and
rational numbers $b_0,c_0$ with denominator $D$ such that
$-M \leq b_0,c_0 \leq M$ and wait for the luck-of-the-draw.

{\bf Removing the Chaff}

Of course $x=1$ is just Gauss. There are also Kummers with $x=-1$ and
their associates via the transformation rules that have to be removed.

{\bf The Wheat}

The output with $M=4$ and $D \leq 4$ and $d=6$ was rather numerous,
but most of them turned out to be special cases of the
following two {\it infinite} families of strange evaluations.

{\bf Theorem 1.} For any non-negative integer $r$, we have
$$
F(-2n,b,-2n+2r-b,-1)=
{{(\khetzi)_n (b+1-r)_n} \over { (b/2+1-r)_n (b/2+\khetzi-r)_n}}
\cdot \left (
\sum_{i=0}^{r-1}  
{  { 2^{2i} i! {{r+i-1} \choose {2i}}  }
\over {(b-r+1)_i} } \cdot {{n} \choose {i}} \right )  \quad .
\eqno(PerturbedKummer)
$$
{\bf Note.} There is an analogous statement with $2r$ replaced by
$2r+1$ that we omit. There are also analogous statements
for $F(-2n,b+r,-2n-b,-1)$ that we also omit).

{\bf Three Sketches of Three Proofs}

{\bf First Proof.} (Shalosh B. Ekhad) For each {\it specific} $r$ this is immediately
doable by the Zeilberger algorithm, but even the general
case is doable that way by dividing both sides of
$(PerturbedKummer)$ by
$$
{{(\khetzi)_n (b+1-r)_n} \over { (b/2+1-r)_n (b/2+\khetzi-r)_n}} \quad,
$$
leaving a polynomial
on the right side. Now apply the Zeilberger algorithm
to this new left side,
getting a certain second-order linear recurrence operator 
that annihilates
this new left side (leaving $r$ symbolic),
and then verifying that
the polynomial on the right side indeeed annihilates it
(and checking the trivial initial conditions).

{\bf Second Proof.} (Dennis Stanton) Use an appropriate specialization
of Eq. (1), Sec. 4.7 of [B], and iterate.

{\bf Third Proof.} (Christian Krattenthaler)
Using his versatile package HYP ([K]), Krattenthaler
used contiguous relations. See his write-up

{\eighttt
http://www.math.rutgers.edu/\Tilde zeilberg/mamarim/mamarimhtml/ck.pdf}, 

that he kindly allowed us to post.

{\bf Conjecture 1.} For any integers $i$ and $j$,
$$
F(-2n,-\khetzi+i,-3n-\khetzi+j,-3) \quad
$$
is evaluable in closed-form. 

{\bf Comment.} We can easily find the explicit forms for
every specific $i$ and $j$. The list of all those
exact evaluations for $-5 \leq i,j \leq 5$ can be found in:

{\eighttt http://www.math.rutgers.edu/\Tilde zeilberg/tokhniot/findhg/oApaZclosedForm} $\quad$ .

However, we were unable to find a
{\it uniform} expression in terms of $i$ and $j$, like in 
Theorem 1. It is very possible that it can be proved
(perhaps without being able to write it uniformly), by
using contiguous relations, as done in Krattenthaler's proof above.

{\bf The Remaining Strange Identities}

After removing all the identities covered by Theorem 1 and
Conjecture 1, 19 (inequivalent) strange hypergeometric 
identities remained. Most of them are already in [E], so
we won't list them here, but refer the reader to {\tt PreComputed21();}
in our package {\tt BruteTwoFone}, and, in a more human-readable form,
complete with the evaluations to:

{\eighttt http://www.math.rutgers.edu/\Tilde zeilberg/tokhniot/findhg/oSefer21} $\quad$ .

{\bf Acknowledgement}

This paper was inspired by an intriguing question posed by
John Greene (to find a one-parameter family generalizing
the exact evaluation, in terms of the Gamma function evaluated 
at rational arguments, of $F({{1} \over {4}}, {{3} \over {4}},1, - {{1} \over {63}})$ ) that
unfortunately we were unable to answer. We also
thank Dennis Stanton and Christian Krattenthaler for
supplying their proofs of Theorem 1, and last but
not least, we would like to thank Shalosh B. Ekhad for
stimulating discussions and insight.

\vfill\eject

{\bf References}

[AAR] G.E. Andrews, R. Askey, and R. Roy, {\it ``Special Functions''},
Cambridge Univ. Press, 1999.

[B] W.N. Bailey, {\it ``Generalized Hypergeometric Series''}, Cambridge University
Press, 1935. Reprinted by Hafner Pub. Co., New York, 1972.

[E] S. B. Ekhad, {\it Forty  ``Strange'' computer-discovered
[and computer-proved (of course!)] hypergeometric series evaluations},
The Personal Journal of Ekhad and Zeilberger, \hfill\break
{\eighttt
http://www.math.rutgers.edu/\Tilde zeilberg/pj.html}, Oct. 12, 2004.

[K] C. Krattenthaler, {\tt HYP}, A Mathematica package 
for the manipulation and identification of binomial and hypergeometric series and 
identities available from \hfill\break
{\tt http://www.mat.univie.ac.at/\Tilde kratt/hyp\_$\,$hypq/hyp.html}

[GS] I. Gessel and D. Stanton, 
{\it Strange evaluations of hypergeometric series},
SIAM J. Math. Anal. {\bf 13}(1982), 295-308.

[MZ] M. Mohammed and D. Zeilberger,
{\it Sharp upper bounds for the orders outputted by the Zeilberger and
q-Zeilberger algorithms}, J. Symbolic Computation {\bf 39} (2005), 201-207.
[available on-line from the authors' websites].

[Z] D. Zeilberger, {\it DECONSTRUCTING the ZEILBERGER algorithm},
J. of Difference Equations and Applications {\bf 11} (2005), 851-856.
[available on-line from the author's website].

\bigskip

{\bf Moa Apagodu}, 
Department of Mathematics, Virginia Commonwealth 
University, Richmond, VA 23284, USA.
\hfill\break
Email: {\eighttt mapagodu  at  vcu dot edu} .
Website: {\eighttt http://www.people.vcu.edu/\Tilde mapagodu/} .
\medskip
{\bf Doron Zeilberger}, 
Department of Mathematics, Rutgers University (New Brunswick),
Hill Center-Busch Campus, 110 Frelinghuysen Rd., Piscataway,
NJ 08854-8019, USA.
\hfill\break
Email: {\eighttt zeilberg at math dot rutgers dot edu} .
Website: {\eighttt http://www.math.rutgers.edu/\Tilde zeilberg/} .
\bigskip
First Written: Feb. 22, 2008

This version:  Feb. 22, 2008
\end